\documentclass[10pt]{amsart}
\usepackage{latexsym, amscd, amsmath, verbatim, amssymb, graphics}
\newtheorem{theorem}{Theorem}[section]
\newtheorem{corollary}[theorem]{Corollary}
\newtheorem{proposition}[theorem]{Proposition}
\newtheorem{lemma}[theorem]{Lemma}
\theoremstyle{remark}
\newtheorem{remark}[theorem]{Remark}
\theoremstyle{definition}
\newtheorem{definition}[theorem]{Definition}

\numberwithin{equation}{theorem}

\newcommand{\Ht}{\operatorname{ht}}

\renewcommand{\hom}{\operatorname{Hom}}
\newcommand{\tor}{\operatorname{Tor}}
\newcommand{\ext}{\operatorname{Ext}}

\newcommand{\depth}{\operatorname{depth}}

\newcommand{\reg}{\operatorname{reg}}

\begin{document}
\title{Bounds on the Castelnuovo-Mumford regularity of tensor products}
\author{Giulio Caviglia}

\begin{abstract}
In this paper we show how, given a complex of graded modules and
knowing some partial Castelnuovo-Mumford regularities for all the
modules in the complex and for all the positive homologies, it is possible
to get a bound on the regularity of the zero homology. We use this to
prove that if $\dim \tor_1^R(M,N)\leq1$ then $\reg(M\otimes N)\leq \reg(
M)+\reg(N)$, generalizing results of Chandler, Conca
and Herzog, and Sidman. Finally we give a description of the regularity of a
module in terms of the postulation numbers of filter regular hyperplane
restrictions.
\end{abstract}
\address{Mathematics Department\\
         UC Berkeley \\
         Berkeley, CA, USA}
\email{caviglia@math.berkeley.edu}
\thanks{The author was partially supported by the ``Istituto Nazionale
di Alta Matematica Francesco Severi'', Rome.} \subjclass{Primary
13D45, 13D02} \keywords{Castelnuovo-Mumford regularity,
postulation number, filter-regular sequence}
 \maketitle

\section{introduction}
Let $R=K[X_1,\dots,X_n]$ be a polynomial ring over a field $K$, $M$ a
finitely generated graded $R$-module and let $I\subset R$ be an ideal. Recently some work has
been done to study when the Castelnuovo-Mumford regularity of $I^r$
can be bounded by $r$ times the regularity of $I,$ and more generally
when the regularity of $IM$ can be bounded by the sum of the regularity of
$I$ and $M.$ This is not always the case;  see the papers
of Sturmfels \cite{St}, and Conca, Herzog \cite{CH} for counterexamples.
On the other hand, under the hypothesis that  $\dim(R/I)\leq 1$, Chandler \cite{Ch} and Geramita, Gimigliano and Pitteloud
\cite{GGP} showed that  $\reg(I^r)\leq
r \reg(I).$ In a recent paper Conca and Herzog \cite{CH} proved,
using similar methods to the one in \cite{Ch} that under the same assumption
(i.e. $\dim(R/I)\leq 1$), $\reg(IM)\leq \reg(I)+\reg(M).$ An extension
of the latter was recently done by Sidman
\cite{S} who showed that if two ideals of $R$, say $I$ and $J$, define
schemes whose intersection is a finite set of points then
$\reg(IJ)\leq \reg(I)+\reg(J).$ She deduced this theorem from a result in
the same paper \cite{S} which bounded the regularity of a tensor product of
sheaves.

In the first section of this paper we show how the same
technique as in \cite{S} can be applied
to prove a stronger statement, i.e that given $M$ and $N$ graded $R$-modules such that $\dim \tor_1^R(M,N)\leq 1,$ then $\reg(M\otimes N)\leq
\reg(M)+\reg(N)$. It is easy to see that this implies all the
other results  mentioned above. This theorem has been recently applied by
Daniel Giaimo \cite{Gi} to prove the Eisenbud-Goto regularity conjecture for
connected absolutely  reduced  curves.

In section 2 we deduce from a formula of Serre that the
Castelnuovo-Mumford regularity can be described in terms of the postulation numbers of filter regular hyperplane
restrictions, where the postulation number $\alpha(M)$ of a module $M$
is defined as the largest nonnegative integer for which the Hilbert
function of $M$ is not equal to the corresponding Hilbert polynomial.
More precisely we show that
given a finitely generated graded $R$-module $M$ with $\dim(M)=d$ we have
\[
\reg(M)=\max_{i \in \{0,\dots,d \}}\{\alpha(M/(l_1,\dots,l_i)M)-\sum_{j=1}^i(D_j-1)\}
\]
where $l_1,\dots, l_d$ is a filter regular sequence on $M$ of degrees
$D_1,\dots, D_d.$

We thank the referee for pointing our attention to the first
section of a paper of Malgrange \cite{Mal}. Looking at another
definition of regularity that goes back to E. Cartan (1901), for
which Janet (1927) proved that $m$-regularity implies
$(m+1)$-regularity, the above theorem, at least for linear forms,
follows from the equivalence to today's definition (attributed by
Malgrange to
 Quilllen, Serre and Mumford in the 1960's). Moreover the extension to higher degree
 forms may also be done along the lines of \cite{Mal}.

\section{Generalities}
From now on by  $R$ we denote the polynomial ring $K[X_1,\dots ,X_n]$ and by $R_+$
its homogeneous maximal ideal.
Given a graded $R$-module $M=\bigoplus_{i\in \mathbb Z} M_i$ we will denote $\bigoplus_{i\in \mathbb Z, i> a }M_i$ by $M_{>a}.$
We want to define some partial Castelnuovo-Mumford regularities for $M$ with respect to a set of 
indices $\mathcal X\subseteq \mathbb {N} $ as follows.

\begin{definition} Given a set of indices $\mathcal X\subseteq
\mathbb {N}$ and a finitely generated graded $R$-module $M$, we
say that $M$ is \emph{$m$-regular with respect to ${\mathcal X}$}
$($i.e.  \emph{$m$-$\reg^{\mathcal X}$$)$} if  we have
$H_{R_+}^i(M)_{>m-i}=0$ for all $i\in \mathcal X$. The
\emph{regularity} of $M$ \emph{with respect to ${\mathcal X}$} is
defined to be the minimum of all the $m$ for which $M$ is
$m$-$\reg^{\mathcal X}.$
\end {definition}
\begin{remark} \label {Grothendieck}
We should observe that, from the Grothendieck vanishing theorem,
all the local cohomology modules are zero for indexes bigger than
$n.$ Note also that when $\mathcal X=\{0,\dots,n\},$ the
$m$-$\reg^{\mathcal X}$ agrees with $m$-regularity in the sense of
Castelnuovo-Mumford.

%therefore it makes sense to use the following notation: given $a\in
%\mathbb Z$ we set $\mathcal X +a$ to  be $\{i+a | i\in \mathcal X
%\}\cap \mathbb \{0,\dots,n\}.$
\end{remark}
The next lemma describes the behavior of the regularity with
respect to ${\mathcal X}$ for exact sequences. We will use the
following notation: given $a\in \mathbb Z$ we set $\mathcal X +a$
to  be $\{i+a | i\in \mathcal X \}\cap \mathbb {N}.$

\begin{lemma}\label{lemma1} Given a short exact sequence of finitely
generated graded $R$-modules,
\[
\begin{CD}
0 @>>> M @>>> N @>>> P @>>> 0,
\end{CD}
\]
we have:
\begin{enumerate}
\item If $M$ and $P$ are $m$-$\reg^{\mathcal X}$ so is $N$. \item
If $N$ is $m$-$\reg^{\mathcal X}$  and $P$ is
$(m-1)$-$\reg^{\mathcal X-1}$,  then $M$ is $m$-$\reg^{\mathcal
X}$. \item If $M$ is $(m+1)$-$\reg^{\mathcal X+1}$ and $N$ is
$m$-$\reg^{\mathcal X}$ , then $P$ is $m$-$\reg^{\mathcal X}$.
\end{enumerate}
\end{lemma}
\begin{proof}
The result follows from the long exact sequence of local
cohomology modules.
\end{proof}
\section{regularity of tensor products and $\hom$ of modules}

The following lemma was inspired by Lemma 1.4 in \cite{S}.

\begin{lemma}\label{lemmacom}Let $\mathbf C$ be a complex of finitely
generated graded $R$-modules
\[
\begin{CD}
\mathbf C:\; 0 @>>> C_{n} @>>> C_{n-1} @>>> \dots @>>> C_0 @>>> 0.
\end{CD}
\]
%\begin{itemize}
%\item
%$\bullet$
If $C_i$ is $(m+i)$-$\reg^{\mathcal X+i}$ for all $i> 0$ and
the $i^{\text{th}}$ homology $H_i(\mathbf C)$ is $(m+i+1)$-$\reg^{\mathcal X+i+1}$ for all $i>0$ then:
\begin{enumerate}
\item[(1)] The $i^{\text {th}}$ boundary $B_i$ is
$(m+i+1)$-$\reg^{\mathcal X+i+1}$ for all $i\geq 0.$ \item[(2)] If
$C_0$ is $m$-$\reg^{\mathcal X}$, then so is $H_0(\mathbf C).$
\end{enumerate}
%\item
%$\bullet$
If $C_{n-i}$ is $(m-i)$-$\reg^{\mathcal X-i}$ for all $i \geq 0$
and the $(n-i)^{\text{th}}$ homology $H_{n-i}(\mathbf C)$ is
$(m-i-1)$-$\reg^{\mathcal X-i-1}$ for all $i>0$, then:
\begin{enumerate}
\item[(1$^\prime$)] The $(n-i)^{\text {th}}$ cycles $Z_{n-i}$ are
$(m-i)$-$\reg^{\mathcal X-i}$ for all $i\geq 0.$ \item[(2$^
\prime$)] In particular $H_n(\mathbf C)$ is $m$-$\reg^{\mathcal
X}.$
\end{enumerate}
%\end{itemize}
\end{lemma}
\begin{proof}
First we prove (1). Note that when $i=n$ the boundary $B_i=B_n=0$ is trivially $(m+i+1)$-$\reg^{\mathcal X+i+1}$.
We can therefore do a reverse induction on $i.$
Consider the following diagram with exact rows and column:
\[
\begin{CD}
@. @. 0 @. @. @. \\
@. @. @VVV @. @. @.\\
0 @>>> B_i @>>> Z_i @>>> H_i( \mathbf C) @>>> 0 @.\\
@. @. @VVV @. @. \\
\dots @>>> C_{i+1} @>>> C_i @>>> C_{i-1} @>>> \dots @.\\
@. @. @VVV @. @. @. \\
@. 0 @>>> B_{i-1} @>>> Z_{i-1} @>>> H_{i-1}(\mathbf C) @>>> 0. \\
@. @. @VVV @. @. @.\\
@. @. 0 @. @. @. @. \\
\end{CD}
\]
By induction we know that
$B_i$ is $(m+i+1)$-$\reg^{\mathcal X+i+1}$ and by assumption
$H_i(\mathbf C)$ is $(m+i+1)$-$\reg^{\mathcal X+i+1}$ so, applying
Lemma \ref{lemma1} to the top exact row in the diagram above, we deduce that
$Z_i$ is $(m+i+1)$-$\reg^{\mathcal X+i+1}.$ Now, since  $C_i$ is
$(m+i)$-$\reg^{\mathcal X+i}$, applying Lemma \ref{lemma1} to the
exact column of the diagram we obtain that $B_{i-1}$ is
$(m+i)$-$\reg^{\mathcal X+i};$ this completes the induction and proves
(1).

We now prove (2).
Consider the exact sequence
\[
\begin{CD}
0 @>>> B_0 @>>> C_0 @>>> H_0 @>>> 0.
\end{CD}
\]
By part (1) we know that $B_0$ is  $(m+1)$-$\reg^{\mathcal X+1}$
and by assumption $C_0$ is $m$-$\reg^{\mathcal X}.$ Therefore from
Lemma \ref{lemma1} it follows that $H_0$ is  $m$-$\reg^{\mathcal
X}.$

The proof of (1$^\prime$) and (2$^\prime$) follow similar lines.
Note that since $Z_n \cong H_n(\mathbf C)$, it is sufficient to
prove (1$^\prime$). Moreover $Z_0=C_0$ is $(m-n)$-$\reg^{\mathcal
X-n}$; we can therefore do a reverse induction on $i$. Apply Lemma
\ref{lemma1}(2) to the last row in the diagram to get $B_{n-i}$ is
$(m-i)$-$\reg^{\mathcal X -i}$ and then apply Lemma \ref{lemma1}
(2) to the exact column to get $Z_{n-i+1}$ is
$(m-i+1)$-$\reg^{\mathcal X-i+1}.$ This completes the induction.
\end{proof}
\subsection {Bounds on the regularity of the tensor product}
$\;$

\vspace{7pt}An easy corollary of Lemma \ref{lemmacom}(2) is the
following.
\begin{theorem}\label{tensor}  Let $M$ and $N$ be finitely generated
graded $R$-modules such that $\mathcal X=\{a,\dots,n\}$, for
$\;a\geq 0$, $M$ is $m$-regular $($i.e.
$m$-$\reg^{\{0,\dots,n\}}$$)$, $N$ is $s$-$\reg^{\mathcal X}$ and
$\tor_i^R(M,N)$ is $(m+s+i+1)$-$\reg^{\mathcal X+i+1}$ for all
$i>0$. Then $M\otimes_R N$ is $(m+s)$-$\reg^{\mathcal X}$.
\end{theorem}
\begin{proof}
Take a minimal graded free resolution $\mathbb F: \dots
\rightarrow\; F_i\; \rightarrow \; \dots \rightarrow \;F_0$ of
$M.$ Note that since $M$ is $m$-regular, the lowest possible shift
appearing in $F_i$ is $-m-i.$ Hence $F_i \otimes N$ is
$(m+s+i)$-$\reg^{\mathcal X}$ and so in particular it is
$(m+s+i)$-$\reg^{\mathcal X+i}.$ The homologies of the complex
$\mathbb F \otimes_R N$ are $\tor_i^R(M,N)$, and by assumption
they are $(m+s+i+1)$-$\reg^{\mathbb X+i+1}$, for $i>0$. The
conclusion follows from Lemma \ref{lemmacom}(2) applied to
$\mathbb F \otimes N$ after noticing that $H_0(\mathbb F \otimes
N)$ is $M\otimes N.$
\end{proof}

\begin{remark}
Note that the condition, ``$\tor_i^R(M,N)$ is
$(m+s+i+1)$-$\reg^{\mathcal X+i+1}$'', of Theorem \ref{tensor} is
clearly satisfied when the Krull dimension of $ \tor_i^R(M,N)$ is less
than or equal to the minimum of $\mathcal X+i$ (since the relevant
local cohomology modules are zero for reasons of dimension).
\end{remark}
Setting $\mathcal X = \{0,\dots,n\}$ and noticing that by rigidity
of $\tor$ (see \cite{A} Theorem 2.1) $\dim \tor_1^R(M,N)\leq 1$ is
equivalent to $\dim \tor_i^R(M,N)\leq 1$ for all $i\geq 1$, we
have the following corollary.

\begin{corollary} \label{regtens}
Let $M$ be an $m$-regular finitely generated graded $R$-module and $N$ be
an $n$-regular finitely generated graded $R$-module such that $\dim \tor_1^R(M,N)\leq 1.$ Then $M\otimes N$ is $(m+n)$-regular.
\end{corollary}

From Corollary \ref{regtens} we can deduce:
\begin{theorem}\label{IM}
Let $I\subseteq R$ be a homogeneous ideal and $M$ a finitely
generated graded $R$-module such that the dimension of
$\tor_1^R(M,R/I)$ is less than or equal to 1. Then $\reg(IM)\leq
\reg(I) + \reg(M).$
\end{theorem}
\begin{proof}
First note that unless $I$ is the whole
ring (in which case the result is obvious), we can assume that
$\reg(I)>0.$ From the exact sequence
\[
\begin{CD}
0 @>>> I @>>> R @>>> R/I @>>> 0
\end{CD}
\]
we get $\reg(R/I)=\reg(I) -1.$
 By Corollary \ref{regtens} $
 \reg(M/IM)=\reg(M \otimes _R R/I)\leq \reg(M) +\reg(I) -1.$
On the other hand applying Lemma \ref{lemma1}(2) to the
 exact sequence
\[
\begin{CD}
0 @>>> IM @>>> M @>>> M/IM @>>> 0,
\end{CD}
\]
we obtain $\reg(IM)\leq \max\{\reg(M),\reg(M/IM)+1\}$ which is less
than or
equal to
\[
\max\{\reg(M),\reg(M)+\reg(I)-1+1\}\leq \reg(M) +\reg(I).
\]
\end{proof}
Theorem \ref{IM} implies the following.

\begin{theorem}[Conca, Herzog, Theorem 2.5 of \cite{CH}]
Let $I\subset R$ be a homogeneous ideal with $\dim R/I\leq 1$ and
let $M$ be a finitely generated graded $R$-module. Then
$\reg(IM)\leq$ $\reg(I)+\reg(M).$
\end{theorem}
\begin{theorem}[Sidman Theorem 1.8 of \cite{S}] Let $I,J$ be
homogeneous ideals of $R$ such that the dimension of $R/(I+J)$ is
less than or equal to 1. Then $\reg(IJ)\leq$ $\reg(I)$+$\reg(J).$
\end{theorem}

\subsection{Bounding the Castelnuovo-Mumford regularity of
$\hom_R(M,N)$}
$\;$

\vspace{7pt} Similar reasoning as in Theorem \ref{tensor} can be
used to prove a bound for the regularity of $\hom_R(M,N),$ where
$M$ and $N$ are finitely generated graded $R$-modules. In this
context the dimensional condition required of $\tor_1^R(M,N)$ has
an analogue in certain conditions on the depth of $\ext^i_R(M,N).$

We prove the following:

\begin{theorem}\label{hom}
Let $M$ and $N$ be finitely generated graded $R$-modules. Let $m$ be
the lowest degree of a homogeneous minimal generator for $M,$ and let $\mathcal
X=\{0,\dots,a\}, \;a\leq n$ be a set of indices. If $N$ is $s$-$\reg^{\mathcal X}$ and $\ext_R^i(M,N)$ is $(s-m-i-1)$-$\reg^{\mathcal X-i-1}$ for all $i>0$, then $\hom_R(M, N)$ is $(s-m)$-$\reg^{\mathcal X}.$
\end{theorem}
\begin{proof} Take a minimal graded free resolution $\mathbb F: \dots \rightarrow\;
F_i\; \rightarrow \; \dots \rightarrow \;F_0$ of $M.$ Note that,
since the lowest degree of a homogeneous minimal generator for $M$
is $m,$ the biggest possible shift appearing in $F_i$ is less than
or equal to $-m-i.$ Hence $\hom_R(F_i , N)$ is
$(s-m-i)$-$\reg^{\mathcal X},$ so in particular it is
$(s-m-i)$-$\reg^{\mathcal X-i}.$ The homologies of the complex
$\hom_R(\mathbb F, N)$ are $\ext_R^i(M,N)$, and by assumption they
are $(s-m-i-1)$-$\reg^{\mathbb X-i-1}$ for all $i>0$. The
conclusion follow from Lemma \ref{lemmacom}(2$^\prime$) applied to
$\hom_R(\mathbb F, N)$ after noticing that $H_n(\hom_R(\mathbb
F,N))$ is $\hom_R(M,N).$
\end{proof}
\begin{remark} The condition: ``$\ext_R^i(M,N)$ is
$(s-m-i-1)$-$\reg^{\mathcal X-i-1}$ for all $i>0$'' of Theorem
\ref{hom} is obtained for example when $\depth \ext_R^i(M,N)$ is
greater than or equal to $n-i$ for all $i>0,$ because in this case
$H_{R_+}^j(\ext_R^i(M,N))=0$ for $j<n-i-1.$  On the other hand,
since for any prime ideal $P$ of $\Ht P<i,$ $\ext_R^i(M,N)_P=0$,
we have $\dim \ext_R^i(M,N)\leq n-i.$ Therefore $\depth
\ext_R^i(M,N)\geq n-i$ if and only if $\ext_R^i(M,N)$ is
Cohen-Macaulay.

Hence we have the following result that is analogous to Theorem
\ref{regtens}.
\end{remark}
\begin{theorem}\label{hom1}
Let $M$ be a finitely generated graded $R$-module with $m$ the
lowest degree of a homogeneous minimal generator of $M,$ and let
$N$ be a finitely generated graded $R$-module such that
$\ext_R^i(M,N)$ is Cohen-Macaulay for all $i>0.$ Then $\reg
(\hom_R(M, N))\leq \reg (N)-m.$
\end{theorem}

\section{The Castelnuovo-Mumford regularity in terms of postulation
numbers of some  hyperplane sections} Let $R$ be
$K[X_1,\dots,X_n]$ with the standard grading and let
$M=\bigoplus_{i \in \mathbb Z} M_i$ be a finitely generated graded
$R$-module of Krull dimension $d$.  In this section we prove how
the Castelnuovo-Mumford regularity of $M$ can be obtained as the
maximum of all the postulation numbers of $d$ filter regular
hyperplane sections. In the following we will denote by $H_M(i)$
the value at $i$ of the Hilbert function of $M$ (i.e.
$H_M(i)=\dim_K M_i$), and with $P_M(i)$ the corresponding Hilbert
polynomial. It is well known that $P_M(i)$ agrees with $H_M(i)$
for $i\gg 0.$ We recall also that, by a theorem of Hilbert, the
Hilbert series (i.e. the formal series defined as $\sum_{i\in
\mathbb Z}H_M(i)Z^i$) has a rational expression
$\frac{h(Z)}{(1-Z)^d},$ where $h(Z)\in \mathbb{Z}[Z,1/Z].$ When a
graded $R$-module $M$ has dimension $0$, we will denote by $\max
M$ the degree of its highest nonzero graded component.

\begin{definition}
Let $M$ be a finitely generated graded $R$-module with Hilbert
series $\frac{h(Z)}{(1-Z)^d}.$ Let $h(Z)=\sum_{i=a}^bc_iZ^i$ with $c_b\not =0.$ We set the \emph {postulation number} of $M$ to be $\alpha(M)=b-d.$
\end{definition}

\begin{remark}It is a well-known fact that the postulation number of
$M$ is equal to the highest degree $i$ for which the Hilbert
function differs from the Hilbert polynomial (i.e.
$H_{M}(i)-P_{M}(i)\not= 0$). For a proof see for example
Proposition 4.1.12 in \cite{BH}. The following formula of Serre
(see \cite {BH}, Theorem 4.4.3, for a proof)
\begin{equation}
H_M(i)-P_M(i)=\sum_{j=0}^d(-1)^j\dim_K H_{R_+}^j(M)_i \mathrm{\;\;for\;all
\;}i\in \mathbb Z, \label{Serre}
\end{equation}
shows how the postulation number of $M$ can be defined in terms of the local
cohomology modules $H_{R_+}^i(M).$
\end{remark}

\begin{definition} A homogeneous element $l\in R$ of degree $D$ is
\emph{filter regular} on a graded $R$-module $M$ if the
multiplication map $l:M_{i-D} \rightarrow M_i$ is injective for
all $i\gg 0$. A sequence $l_1,\dots,l_m$ of homogeneous elements
of $R$ is called a \emph{filter regular sequence} on $M$ if $l_i$
is filter regular on $M/(l_1,\dots,l_i-1)M$ for $i=1,\dots,m.$
\end{definition}
\begin{remark}
Since $H_{R_+}^0(M)=\{u\in M\;|\; (R_+)^ku=0$ for some $k\},$ then
$l$ is filter regular on $M$ if and only if $l$ is a
nonzerodivisor on $M/H_{R_+}^0(M).$
\end{remark}
The proposition below will be useful in the sequel. The proof uses a slight modification
of the argument in \cite{M} (p. 101-102). Moreover,  some of its  corollaries (Corollary \ref{b},
 \ Corollary \ref{Conca} and also Theorem \ref{last}) can also be obtained using the same
 method as in \cite{M}.
\begin{proposition} \label{almost}
Let $M$ be a finitely generated graded $R$-module and let $l\in R$
be a filter regular element on $M$ of degree $D.$ Then for any set
of indices ${\mathcal X}\subseteq \{0,\dots,n\}$ we have:
\begin{enumerate}
\item
$\reg^{\mathcal X+1}(M)\leq \reg^{\mathcal X\cup ({\mathcal X}+1)}(M/lM)-D+1$
\item
$\reg^{\mathcal X}(M/lM)-D+1 \leq \reg^{{\mathcal X} \cup ({\mathcal X}+1)}(M).
$
\end{enumerate}
\end{proposition}

\begin{proof}
Consider the short exact sequence
\[
\begin{CD}
0 @>>> M/0:_M l(-d) @>\cdot l>> M @>>> M/lM @>>> 0.
\end{CD}
\]
Note that, since $l$ is filter regular on $M$, $H_{R_+}^i(M/0:_M
l(-D))\cong H_{R_+}^i(M)(-D)$ for all $i>0.$
Looking at the long exact sequence of local cohomology modules, we have
\[
\begin{CD}
\dots @>>> H_{R_+}^i(M) @>>> H_{R_+}^i(M/lM) @>>> H_{R_+}^{i+1}(M)(-D)
@>>> \\
@>>>  H_{R_+}^{i+1}(M)
@>>> H_{R_+}^{i+1}(M/lM)          @>>> \dots
\end{CD}
\]
for all $i \geq 0.$

Let $j> \reg^{\mathcal X}(M/lM)-D+1$ and let $i\in {\mathcal X}.$
Consider the exact sequence
 of $K$-vector spaces given by the graded pieces of degree $j-i+D-1$ of the
 previous sequence. Because of the choice of $j$, we have
 $H_{R_+}^i(M/lM)_{j-i+D-1}=H_{R_+}^{i+1}(M/lM)_{j-i+D-1}=0.$
 Therefore,
 $$H_{R_+}^{i+1}(M)(-D)_{j-i+D-1}\cong
 H_{R_+}^{i+1}(M)_{j-i+D-1},$$
 that is
 $$H_{R_+}^{i+1}(M)_{j-i-1}\cong
 H_{R_+}^{i+1}(M)_{j-i+D-1}.$$
 After a simple induction,
 $H_{R_+}^{i+1}(M)_{j-i-1}\cong  H_{R_+}^{i+1}(M)_{j-i+sD-1}$ for any
 $s>0.$ Since  $H_{R_+}^{i+1}(M)$ is Artinian, we obtain that $H_{R_+}^{i+1}(M)_{j-i-1}=0$ for all $i\in {\mathcal X},$ which implies part (1).

We now prove part (2). Take $j>\reg^{{\mathcal X} \cup ({\mathcal
X}+1)}(M)+D-1$ and $i\in {\mathcal X}.$ From the choice of $j$, we
have $H_{R_+}^i(M)_{j-i}=H_{R_+}^{i+1}(M)(-D)_{j-i}=0$ for any
$i\in {\mathcal X}.$ In particular looking at the $(j-i)^{\mathrm
{th}}$ graded component of the long exact sequence of local
cohomology modules, we get $H_{R_+}^i(M/lM)_{j-i}=0$ for all $i\in
{\mathcal X}$, which implies part (2).
\end{proof}

Proposition \ref{almost} has the following corollaries:
\begin{corollary}\label{b}
Given a finitely generated graded $R$-module $M$ and a filter
regular element $l$ of degree $D,$ we have
\[
\reg(M/H_{R_+}^0(M))\leq \reg(M/lM)-D+1.
\]
\end{corollary}
\begin{proof}
Set ${\mathcal X}=\{0,\dots,n\}$ and note that
$\reg(M/H_{R_+}^0(M))=\reg^{\mathcal X+1}(M).$ The conclusion
follows from Proposition \ref{almost}(1).
\end{proof}
\begin{corollary}[\cite{CH} Proposition 1.2]\label{Conca}
Given a finitely generated graded $R$-module $M$ and a filter
regular element $l$ of degree $D$ we have
\[
\reg(M)=\max\{\max H_{R_+}^0(M), \reg(M/lM)-D+1\}.
\]
\end{corollary}
\begin{proof}
Take ${\mathcal X}=\{0,\dots,n\},$ and note that
$\reg(M)=\max\{\reg^{\{0\}}(M),\reg^{{\mathcal X}+1}(M)\}.$
Clearly $\reg^{\{0\}}(M)=\max H_{R_+}^0(M).$ From Proposition
\ref{almost} (1) we have $\reg^{{\mathcal X}+1}(M)\leq
\reg^{\mathcal X\cup ({\mathcal X}+1)}(M/lM)-D+1=\reg(M/lM)-D+1.$
Thus we get $\reg(M)\leq \max\{\max H_{R_+}^0(M),
\reg(M/lM)-D+1\}.$ On the other hand, $\max
H_{R_+}^0(M)\leq\reg(M)$ and, by Proposition \ref{almost}(1), we
have  $\reg^{\mathcal X}(M/lM)-D+1\leq\reg^{\mathcal X\cup
({\mathcal X}+1)}(M)$ $=\reg(M).$
\end{proof}

\begin{theorem}\label {postulation}
Let $M$ be a finitely generated graded $R$-module with $\dim(M)=d.$ Then
\[
\reg (M)=\max_{i \in \{0,\dots,d
\}}\{\alpha(M/(l_1,\dots,l_i)M)-\sum_{j=1}^i(D_j-1)\}
\]
where $l_1,\dots ,l_d$ is a filter regular sequence of degrees $D_1,\dots, D_d.$
\end{theorem}
\begin{proof}

By definition, given any finitely generated graded $R$-module $N$ and
any $i>\reg(N),$ we have $H_{R_+}^j(N)_{i-j}=0.$ In particular
$H_{R_+}^j(N)_{i}=0,$ hence from (\ref{Serre}) it is clear that
$\reg(N)\geq \alpha(N)$ for every $N.$

By Corollary \ref{Conca}, $\reg(M)\geq \reg(M/lM)-\deg l +1$ for any filter
 regular element $l,$ so in particular we have
\[
\reg (M)\geq\max_{i \in \{0,\dots,d
\}}\{\alpha(M/(l_1,\dots,l_i)M)-\sum_{j=1}^i(D_j-1)\}.
\]
We need to prove the reverse inequality. We do an induction on the
dimension of $M.$ If $\dim M=0,$ then $\reg (M)=\max H_{R_+}^0(M)$
which equals to $\alpha(M)$, by (\ref{Serre}). Assume $d>0.$ By
induction hypothesis we get
\[
\reg (M/l_1M)=\max_{i \in \{1,\dots,d
\}}\{\alpha(M/(l_1,l_2,\dots,l_i)M)-\sum_{j=2}^i(D_j-1)\}.
\]
Consequently setting $a=\max_{i \in \{0,\dots,d
\}}\{\alpha(M/(l_1,\dots,l_i)M)-\sum_{j=1}^i(D_j-1)\}$ we have
\[
\reg (M/l_1M)-D_1+1\leq a.
\]
Thanks to Corollary \ref{Conca} we still have to prove that
$\max H_{R_+}^0(M)\leq a.$ Since $H_{R_+}^j(M)\cong
H_{R_+}^j(M/H_{R_+}^0(M))$ for all $j>0$, by  Corollary  \ref{b} we
know that $H_{R_+}^j(M)_{>a-j}=0$ for all $j>0.$ In particular for any $b>a,$ $H_{R_+}^j(M)_b=0$ for all $j>0.$ Hence by (\ref{Serre}) we have $H_{M}(b)-P_{M}(b)=\dim_K H_{R_+}^0(M)_b.$ But $a\geq \alpha(M)$ so $H_{M}(b)-P_{M}(b)=0$ for all $b>a\geq\alpha(M).$ Therefore $\max  H_{R_+}^0(M)\leq a.$
\end{proof}
An interesting corollary of the Theorem \ref{postulation} is the following
\begin{corollary}\label{indip}
Let $M$ be a finitely generated graded $R$-module such that $\dim
M=d,$ and let $l_1,\dots,l_d$ be a filter regular sequence on $M$
of elements of degree $D_1,\dots,D_d.$  Then  the number
\[
\max_{i \in \{0,\dots,d
\}}\{\alpha(M/(l_1,\dots,l_i)M)-\sum_{j=1}^i(D_j-1)\}
\]
is independent of the choice of the filter regular sequence and of its degrees.
\end{corollary}

\begin{remark}
Note that both Theorem \ref{postulation} and Corollary \ref{indip}
lie on the definition of $\alpha(M).$ The number $\alpha(M)$ is the highest integer $i$ for which the function $\phi$ defined as
\[
\phi(i,M_0,M_1,M_2,\dots,M_n):=       \sum_{j=0}^n(-1)^j\dim_K(M_j)_i
\]
is not zero at $(i,H^0_{R_+}(M),H^1_{R_+}(M),\dots,H^n_{R_+}(M)).$
We want to point out that we can substitute for $\phi$ any other
function $\psi$ such that, whenever $(M_j)_{>i-j}=0$ for all $j>0$, we have:
\begin{equation}
\psi(i,M_0,M_1,M_2,\dots,M_n)\not =0 \text{ if and only if } (M_0)_i \not = 0.\label{condition}
\end{equation}
Then instead of $\alpha(M)$ we could use the function $\beta(M)$ defined as:
\[
\sup \{i\;|\;\psi(i,H^0_{R_+}(M),H^1_{R_+}(M),\dots,H^n_{R_+}(M))\not =0\}.
\]
If we set for example $\psi(i,M_0,\dots,M_n)= \dim_K(H^0(M_0)_i),$
we obtain, as an analogue of Theorem \ref{postulation}, the
already known fact:
\begin{theorem}\label{last}
Given a finitely generated module $M$ of dimension $d$ we have
\[
\reg(M)=\max_{i \in \{0,\dots,d
\}}\{sat(M/(l_1,\dots,l_i)M)-\sum_{j=1}^i(D_j-1)\}.
\]
Where $sat(P)$ is defined to be $\max H_{R_+}^0(P).$
\end{theorem}
Note that Theorem \ref{last} can be found in \cite{Gr} (see Theorem
2.30 (5),(6)) under the more
restricted assumptions that the field $K$ has characteristic zero and the
$l_i$'s are generic linear forms. It can also be easily derived from \cite{CH} Proposition 1.2.
\end{remark}

\section*{acknowledgements}
I would like to thank Professor Craig Huneke for many helpful comments concerning this paper.

\end{document}